\documentclass[12pt]{article}

\usepackage[cp1251]{inputenc}
\usepackage[T2A]{fontenc}

\usepackage{amssymb}
\usepackage{amsmath}
\usepackage{amsthm} 
\usepackage{mathtext}
\usepackage{amsfonts}

\newtheorem{theorem}{Theorem}

\theoremstyle{definition}
\newtheorem{remark}{Remark}

\begin{document} 

\title{On solutions of the modified\\ Helmholtz equation in a planar 
domain}
\author{Alexander~I.~Nazarov\footnote{PDMI RAS and St. Petersburg State 
University}
\date{}}
\maketitle

\hfill{\it Dedicated to the memory of Nikolay G. Kuznetsov}

\abstract{A criterion for the convexity of a two-dimensional domain is proved 
in terms of the solution of a boundary value problem for the modified Helmholtz 
equation.}

\section{Introduction}

N.G. Kuznetsov \cite{Kuz1} (see also \cite[Theorem 2.3]{Kuz2}) 
has proved the following statement:
\medskip

\noindent{\it Let $\Omega\subset \mathbb R^2$ be a bounded domain with $\mathcal 
C^2$-boundary. Consider the boundary value problem
\begin{equation}\label{Helm}
 \Delta v_\mu-\mu^2v_\mu=0 \quad \mbox{in} \quad\Omega, \qquad 
v_\mu|_{\partial\Omega}=1.
\end{equation}
If the solutiom of this problem satisfies\footnote{\ By the maximim 
principle (see, e.g., Theorem 2.1 in the survey \cite{AN}, whеre also the 
history of the issue is discussed in detail), the solution of (\ref{Helm}) 
is unique, and
$0<v_\mu<1$ in $\Omega$. }
\begin{equation}\label{ineq}
 |\nabla v_\mu(x)| \leqslant\mu v_\mu(x), \quad x\in \Omega
\end{equation}
for all sufficiently large $\mu>0$, then 
$\Omega$ is convex.
}
\medskip

A natural question arises: is the condition (\ref{ineq}) a criterion for the 
convexity of a domain in $\mathbb R^2$? A positive answer to this question was 
conjectured in \cite{Kuz2}. In this note, this conjecture is proved, even 
without the assumption that $\Omega$ is smooth.

\begin{theorem}
 Let $\Omega\subset \mathbb R^2$ be a boinded Lipschitz domain, and let $v_\mu$ 
be ths solution of the boundary value problem (\ref{Helm}). The inequality
(\ref{ineq}) is satisfied for all sufficiently large $\mu>0$ if and only if 
$\Omega$ is convex. Moreover, in this case the inequality holds for all $\mu>0$. 
\end{theorem}

\section{The proof of ``only if''}

To begin, let us briefly recall the idea of the proof for a smooth domain. We 
put 
$w_\mu(x)=-\mu^{-1}\ln v_\mu(x)$. Then
$$
-\Delta w_\mu=\frac {\Delta v_\mu}{v_\mu} - \frac {|\nabla 
v_\mu|^2}{v_\mu^2}=\mu^2 - \frac {|\nabla 
v_\mu|^2}{v_\mu^2}\geqslant0
$$
for all sufficiently large $\mu$ due to (\ref{ineq}). Thus, $w_\mu$ is 
superharmonic.

By the Varadhan theorem\footnote{\ Just this statement uses the smoothness of
$\partial\Omega$.} \cite[p. 434]{V}, 
$$
w_\mu(x)\rightrightarrows {\rm 
dist}(x,\partial\Omega)\quad \mbox{in}\quad \overline\Omega \quad 
\mbox{as}\quad \mu\to+\infty.
$$
 By direct and inverse theorems on the mean value property for 
(sub/super)harmonic functions, the limit function is also 
superharmonic.\footnote{\ In \cite{Kuz1} the argument here is somewhat more 
complicated.} By the result of \cite{AK}, this implies convexity of 
$\Omega$.\medskip

Now we remove the smoothness condition for $\partial\Omega$. Suppose $\Omega$ is 
not convex. Then there exists a segment with endpoints on $\partial\Omega$ that 
lies outside $\overline{\Omega}$ (except for the endpoints). We may assume that 
this segment lies on the $Ox_1$-axis, and the domain is located in the upper 
half-plane. 

Take a circle of large radius enclosing 
$\overline{\Omega}$ and move it downward until it touches $\partial\Omega$ from 
the inside at some point $x^0$, whose projection onto the $Ox_1$-axis lies on 
this segment. Next, taking a small arc of this circle and smoothly continuing it 
inside $\Omega$, we obtain a domain 
$\Omega'\subset\Omega$ with a smooth 
boundary such that $\partial\Omega'$ touches $\partial\Omega$ at the point $x^0$ 
and has a non-zero (negative) curvature at $x^0$ (we denote it $-\varkappa$).

By \cite[Lemma 14.16]{GT}, the function $d(x)={\rm dist}(x,\partial\Omega')$ 
is smooth in a boundary strip of 
$\Omega'$, and Lemma 14.17 
\cite{GT} gives $\Delta d(x^0)=\varkappa$. So, $\Delta d>0$ in some 
neighborhood of  $x^0$ in $\Omega'$.

Consider the solutiom $v'_\mu$ of the problem (\ref{Helm}) in $\Omega'$. If 
$v'_\mu$ satisfies the inequality (\ref{ineq}) in a neighborhood of the point 
$x^0$ for all sufficiently large $\mu>0$, then the function
$w'_\mu(x)=-\mu^{-1}\ln 
v'_\mu(x)$ is superharmonic in this neighborhood. Then, by Varadhan's theorem,  
$d(x)$ also should be superharmonic, a contradiction. Therefore, 
$$
\partial_{\bf n} v'_\mu(x^0)=-|\nabla v'_\mu(x^0)|<-\mu
$$
(here $\bf n$ is the unit vector of the inward normal to $\partial\Omega'$ at 
the point $x^0$).

Since 
$$
v_\mu|_{\partial\Omega'}\leqslant 1= v'_\mu|_{\partial\Omega'},
$$
and there is no identical equality, the maximum principle gives
$v_\mu<v'_\mu$ in $\Omega'$. But
 $v_\mu(x^0)= v'_\mu(x^0)=1$, therefore
$$
\limsup\limits_{\delta\to +0}\delta^{-1}\big(v_\mu(x^0+\delta{\bf 
n})-v_\mu(x^0)\big)\leqslant \limsup\limits_{\delta\to 
+0}\delta^{-1}\big(v'_\mu(x^0+\delta{\bf 
n})-v'_\mu(x^0)\big)<-\mu,
$$
whence
$$
\liminf\limits_{\delta\to +0}|\nabla v'_\mu(x^0+\delta{\bf n})| \geqslant -
\limsup\limits_{\delta\to +0}\partial_{\bf n} v_\mu(x^0+\delta{\bf 
n})>\mu v_\mu(x^0).
$$
Thus, the inequality (\ref{ineq}) fails in a neighborhood of the point $x^0$, 
and the statement ``only if'' is proved.

\begin{remark}
Note that the Lipschitz condition for 
$\partial\Omega$ was not actually used. A 
priori, the only requirement for the domain is that the problem (\ref{Helm}) is 
solvable.
\end{remark}

\section{The proof of ``if''}

Let $\Omega$ be a convex bounded domain, and let $v_\mu$ be the solution of the 
problem (\ref{Helm}). First, we claim that the inequality (\ref{ineq}) holds on
$\partial\Omega$.

Let $x^0\in\partial\Omega$. We may assume that $x^0=0$, and the domain is 
located in the upper half-plane. We introduce the barrier function
$V_\mu(x)=\exp(-\mu x_2)$.

Since
$$
V_\mu|_{\partial\Omega}\leqslant 1= v_\mu|_{\partial\Omega},
$$
and there is no identical equality, the maximum principle gives
$V_\mu<v_\mu$ in $\Omega$. But $V_\mu(0)= v_\mu(0)=1$, therefore 
$$
|\nabla v_\mu(0)|=-\partial_{x_2}v_\mu(0) \leqslant 
-\partial_{x_2}V_\mu(0)=\mu v_\mu(0),
$$
and the claim follows.

To complete the proof, it suffices to show that all sublevel sets
$$
\Omega_A=\{x\in\Omega:\,v_\mu(x)< A\},\qquad \min\limits_\Omega v_\mu(x)<A<1,
$$ 
are convex. Indeed, in this case, applying the previous argument to the function 
$A^{-1}v_\mu(x)$ in the domain $\Omega_A$, we obtain the desired result.

For the proof, we need the following statement (\cite[Theorem 
4.2]{Ken}):\medskip

\noindent{\it Let $n\geqslant2$, and let 
$\Omega\subset\mathbb R^n$ be a bounded convex domain with $\mathcal 
C^2$-boundary. Let the function 
$u\in \mathcal C(\overline{\Omega})\cap\mathcal C^2(\Omega)$ be a solution of 
the boundary value problem
\begin{equation}\label{Ken}
\Delta u + h(u)=0 \quad \mbox{in} \quad\Omega, \qquad 
u|_{\partial\Omega}=0, \qquad u>0 \quad \mbox{in} \quad\Omega.
\end{equation}
Let the function $h:\mathbb R_+\to\mathbb R_+$ satisfy the following 
assumptions for some $0<\alpha<1$:
\begin{enumerate}
 \item $t^{\alpha-1}h(t)$ is a strictly decreasing function of $t$ variable;
 \item $t^{3-\frac 1\alpha} h(t^{\frac 1\alpha})$ is a concave function of $t$ 
variable. For a twice differentiable $h$ this assumption can be rewritten as 
follows:
$$
(1-2\alpha)(1-3\alpha)h(t)+(5\alpha-1)t h'(t)+t^2h''(t)\leqslant0.
$$
\end{enumerate}
Then the function $u^\alpha$ is concave in $\overline{\Omega}$.
}
\medskip

We temporarily suppose that $\Omega$ is smooth. Consider the function
$u(x)=1-v_\mu(x)$. It solves the problem (\ref{Ken}) with $h(t)=\mu^2(1-t)$. 
This function evidently satisfies\footnote{Formally speaking, the function $h$ 
is positive only on the interval $(0,1)$. But since $u$ takes values only in 
this interval, this does not change the conclusion.} assumptions $1$ and $2$ 
with $\alpha=\frac 12$. Therefore, the function $\sqrt{1-v_\mu}$ is concave, 
and all sets $\Omega_A$ are convex.

It remains to get rid of the smoothness assumption of $\Omega$. But any convex 
domain can be approximated from the inside by smooth convex ones. Since the 
value of $v_\mu$ at any interior point depends continuously on $\partial\Omega$, 
the set $\Omega_A$ for any fixed $A$ is approximated by convex sets and is 
therefore convex itself. So, the statement ``if'' is also proved.

\end{document}